\documentclass[a4paper]{article}
\usepackage[utf8]{inputenc}

\RequirePackage{amsmath}
\RequirePackage{relsize}
\RequirePackage{amsthm}
\RequirePackage{amsfonts}
\RequirePackage{amssymb}
\RequirePackage{url}

\RequirePackage{physics}

\usepackage{physics}
\usepackage{setspace}
\usepackage{graphicx}
\usepackage[font=small,labelfont=bf]{caption}
\usepackage{epstopdf}
\AppendGraphicsExtensions{.pdf}

\usepackage[marginratio=1:1]{geometry}

\newtheorem{remark}{Remark}[section]
\newtheorem*{cremark}{Concluding Remark}
\newtheorem{lemma}{Lemma}[section]

\newtheorem{example}{Example}[section]

\newtheorem{corollary}{Corollary}[section]
\newtheorem{theorem}{Theorem}[section]
\newtheorem{conjecture}{Conjectre}[section]
\newtheorem{proposition}{Proposition}[section]

\author{Alexander Kushkuley \\ kushkuley@gmail.com}
\title{  Absolute values and  tensor powers of irreducible characters      }
\begin{document}
    \maketitle

    \begin{abstract}
        \noindent
        Let  $ \chi $ be a  character of a complex irreducible representation of a finite group $G$.  We present a simple formula for the expectation of the random variable  $(|\chi|/\chi(1))^{t} $ in terms  of character ratios  $  (|\chi(g)|/\chi(1))^{t}, \; g \in G, \; t \geq 0   $. As a follow up we briefly discuss asymptotic properties of the formula and its relation to the subject of growth of dimensions of isotypic components in (virtual)  tensor powers of irreducible representations. Similar type of reasoning can be applied to some questions related  to commuting probability. In particular, we obtain an analogue of Frobenis formula for the probability of an "event" $ |\chi( [x,y]^{-1}g)| = r $
    \end{abstract}
    \onehalfspacing
    \section{Introduction}
    \numberwithin{equation}{section}
    A problem of estimating upper bounds of absolute values of irreducible characters of finite groups (cf. e.g. \cite{GLT1}) is related to a problem of estimating  multiplicities of irreducible components in tensor powers (cf. \cite{COUL}-\cite{NAI}).
    In a  simple  case of representations over the field of complex numbers
    these two questions can be addressed simultaneously.
    It turns out that both estimates depend on expectation of an absolute value of a relevant irreducible character.
    Along these lines, one obtains simple lower and upper bounds on an absolute value of an irreducible character that are better than  "centralizer bound" (mentioned  in e.g. \cite{GLT1}). The same computation "in reverse"  yields some results that resemble (asymptotic) estimates for the total number of irreducible summands in tensor powers  (cf. \cite{He}, \cite{COUL}, \cite{NAI}). A generic "computation" suggested in this paper can be also applied to some questions related to commuting probability (cf. \cite{Howe}, \cite{Serre2}, \cite{prob}). In particular, we  obtain a generalization of the Frobenius character ratio formula (cf. e.g. \cite{Howe}, \cite{Serre2})

    The paper is structured as follows. The introduction contains relevant terminology,
    basic notations, a statement of the main result and some examples.
    Auxiliary lemmas and definitions are collected in section 2.
    A proof of the main result on character value bounds is the subject of section 3 and
     "multiplicities of irreps in tensor powers" are briefly discussed in section 4.
      Some applications to commuting probability  are presented in section 5.

    Let $G$ be a finite non-abelian group and let $ \chi$  be a  character of an irreducible complex   representation
    \begin{align}
        \rho: G \rightarrow \textnormal{GL}(V), \; \dim V > 1   \tag{1.0}
    \end{align}
    The maximal value of the
    central function $|\chi|$ is $ n = \dim V (\equiv \dim \rho \equiv \dim \chi  ) $. Let $ \gamma$ be the second largest value of    $ |\chi| $.
    Set $ K
    = \{ g \in G \; | \;  |\chi(g)| = n \}  $, in other words let  $K$  be  a normal subgroup of all $ G $-elements that act on $V$ by scalar  multiplication (in particular $K \supset \ker \rho $). Suppose that $ \gamma > 0 $ and
    set
    \begin{align}
        C = G \setminus K, \;
        C_0 = \{ g \in C \; | \; |\chi(g)| = \gamma  \}, \; G_0 = |\chi^{-1}(0)| \nonumber
    \end{align}

    \noindent Let
    \begin{align}
        n > \gamma = \gamma_0 > \gamma_1 > \cdots > \gamma_l = 0, \; \gamma \neq 0
    \end{align} be the list of all distinct values of the
    central function $ | \chi | $.
    In accordance with (1.1), set
    \begin{align}
        C_q = \{ g \in C \; | \; |\chi(g)| = \gamma_q \} , \; q = 0, \cdots,  l
    \end{align}
    \noindent In other words, $ K $ and $C_q, \; q = 0, \cdots, l  $ are (all) level sets of the function $ |\chi| $. As a matter of convenience we exclude the "kernel"  $ K $ from the list (1.2) (cf., however, Remarks 2.1 and 5.3)
    \begin{remark}\
        \begin{enumerate}
            \item[(a)]
            $C_q$ is a disjoint union of conjugate classes on which the character   $  |\chi|^2 $
            takes the value $ \gamma_q^{2} $.
            The maximal nontrivial character value $\gamma$ is not necessary non-zero. For example, $ \gamma_0 \equiv \gamma = 0$ for any nonlinear irrep of a non-abelian finite two step nilpotent group (see e.g. \cite{Drinfeld}, Appendix B). We exclude such cases.
            \item[(b)]
            Since $ |\chi(g)| = |\chi(g^{-1})|$, $ C_q $ is a conjugate class only if $ C_q^{-1} = C_q $, i.e. only if $C_q $ is a real conjugacy class (cf. e.g. \cite{Serre})
            \item[(c)]    $ C_0 $ is a  subset of $ G \setminus G_0 \setminus K $ and therefore $ |C_0| \leq |G| - |G_0| -|K|  $. Also, note that $ C_l = G_0 $ and that
            $ \sum_{i=0}^{l-1} |C_i| = |G| - |G_0| -|K| $
        \end{enumerate}
    \end{remark}

    \noindent
    Any finite set $ X $ can be viewed as a probability space with  probability  $\textnormal{Pr}_X(x)$ of $ x \in X $  defined as $ 1/|X| $. An expectation of a  random variable (function) $ \phi : X \rightarrow \mathbb{C}$ will be denoted by  $\mathbb{E}_X(\phi)$. Note, that the standard scalar product
    $$ <\phi, \psi > = \frac{1}{G}\sum_{g \in G} \phi(g)\overline{\psi(g)} $$
    on the space   of complex functions  on $ G $ is an expectation of the function $ \phi \bar{\psi}$. Formally speaking, we have by definition

    \begin{lemma} (cf. e.g. \cite{Serre}).
        \begin{enumerate}
            \item[(i)]
            If $ \phi_1, \phi_2 $ are complex functions on $ G $ then
            $\mathbb{E}_G(\phi_1 \bar{\phi}_2) = \; <\phi_1, \phi_2 > $, in particular, expectation of a character is equal to the dimension of the space of fixed points of corresponding representation.

            \item[(ii)] Suppose that $ \phi = a_1 \chi_1 + \cdots + a_k \chi_s  $ where $ \chi_i, \; i = 1, \cdots s $ are pairwise distinct irreducible characters of $ G $, $ a_i \in \mathbb{C} $. Assume that $ \chi_1$ is the trivial character. Then

            $$ \mathbb{E}_G( \phi ) = a_1 $$
            and
            \begin{align} \mathbb{E}_G( \phi\bar{\phi}) =  \sum_{i=1}^k |a_i|^2    \equiv  \norm{\phi}^2 \nonumber
            \end{align}
        \end{enumerate}
    \end{lemma}
    \noindent The $G$-character
    $    |\chi|^2 \equiv \chi  \bar{\chi} $ is a lift (to $G$) of a character of the factor group $ G/K $.
    Let
    \begin{align}
        \chi_1, \cdots, \chi_k
    \end{align}
    be a list of all pairwise inequivalent
    irreducible complex characters of $ G $ that are trivial on $ K $. Set $ \dim \chi_i = \chi_i(K)  = n_i, \; i=1, \cdots, k $ and assume that $\chi_1$ is the trivial character (hence $ n_1 = 1$).   For any real $ t \geq 0 $ the central function
    $	|\chi|^{t} : G \rightarrow G/K \rightarrow \mathbb{C} $ has  unique decomposition
    \begin{align}
        |\chi|^{t} \equiv (\chi  \bar{\chi})^{t/2} = a_{1,t} \chi_1 + \cdots + a_{k, t} \chi_k
    \end{align}
    where  $ a_{i, t}$  are uniquely defined real numbers (see Lemma 2.4) that are non-negative integers when $t$ is an even integer. In  particular,  $ a_{1, t}$ as the multiplicity  of a trivial character in $ (\chi  \bar{\chi})^{t/2}$  is equal to the dimension of the space of $G$-invariants in  $ (V \otimes V^{*})^{t/2} $ for even integer $t$.
    \begin{remark}
        If  $\rho$ is exact and $G$ does not have a nontrivial center then  $ \chi $ occurs in the list (1.3). Note also that by continuity, $ |\chi|^{0}$ is a characteristic function of the subset $ G \setminus G_0 \subset G  $
    \end{remark}

    \noindent By Lemma 1 (i) $ a_{1,t} $ is an expectation of a non-negative random variable ($|\chi|^t$) and therefore is non-negative real number for any $ t \geq  0 $ (see also Lemma 2.2). Let's agree to denote (non-negative, integer) Kronecker coefficients $ a_{i, 2} $ by $ a_i, \; i = 1, \cdots k $

    \begin{theorem}\
        \begin{enumerate}
            \item[(a)] (cf.  (\cite{GLT1}). For the maximal  nontrivial absolute value $\gamma$ of the character $ \chi$ and any $ t > 0 $
            we have lower and upper bounds
            \begin{align}
                \left(\frac{   a_{1,t}/n^{t} - |K|/|G| }
                {( |G| - |G_0| - |K| )   /|G|} \right)
                ^{\frac{1}{t}}
                \; \leq \;
                \frac{\gamma}{n} \; \leq \;
                \left(\frac{   a_{1,t}/ n^{t} -  |K|/|G| }  {|C_0| /|G|} \right)
                ^{\frac{1}{t}}
            \end{align}
            \noindent or equivalently
            \begin{align}
                \!\!\!\!      \!\!\!\! \!\!\!\!      \!\!\!\! \left(\left(\frac{\gamma}{n} \right)^{t}
                \frac{|G| - |G_0| - |K|}{|G|}  +  \frac{|K|}{|G|} \right)^{1/t}
                \; \geq \;
                \frac{\sqrt[t]{a_{1,t}}}{n} \; \geq \;
                \left(   \left(\frac{\gamma}{n}\right)^{t}  \frac{|C_0|}{|G|} +    \frac{|K|}{|G|}\right)^{1/t}   \tag{1.5'}
            \end{align}
            \item[(b)] Both  sides of (1.5) converge to $ \gamma/n$ when $ t $ goes to infinity
            \item[(c)] (cf. \cite{COUL}-\cite{NAI}). When $t$ approaches infinity, both sides of (1.5') obviously converge to one and, therefore
            $$
            \lim_{ t \rightarrow \infty } (a_{1,t})^{1/t} = n
            $$
            \item[(d)] (cf. \cite{He}).    $ \left(\frac{\gamma}{n}\right)^{t} $ and
            $ \frac {a_{1,t}/n^{t} - |K|/|G|} {|C_0|/|G|}$
            are asymptotically equal  as functions of $ t $
        \end{enumerate}
    \end{theorem}

    \noindent A simple proof of this result will be presented below.
    We would like, however, to point  out that the upper bound in (1.5)  follows directly from Markov inequality (cf. e.g. \cite{B}) as follows
    \begin{align}  \frac{|C_0|} {|C|} \; = \; \textnormal{Pr}_{C}( |\chi| \; \geq \; \gamma ) \; \leq \;  \frac{1} {\gamma^{t} } \mathbb{E}_{C} (|\chi|^{t}) \;
        = \; \frac{  |G| \mathbb{E}_{G} (|\chi|^{t} ) - |K| n^{t}  } {\gamma^{t} |C|  } \nonumber
    \end{align}
    for any $ t > 0  $ where $ \mathbb{E}_{G} (|\chi|^{t} ) = a_{1,t} $ by Lemma 1.1.

    \begin{corollary} $ a_{1,t}/n^{t} \; \geq \; |K|/|G| $ for any $ t \geq 0 $
    \end{corollary}
    \noindent Another corollary (from the proof) of Theorem 1.1 (see section 3) provides upper bounds for all absolute values (1.1) of the character $\chi$
    \begin{corollary}
        \begin{align}
            \gamma_i \; \leq \;
            \left(\frac{   a_{1,t} - n^{t} |K|/|G| }  {
                ( |C_0| + \cdots + |C_i|  )/|G|}  \right)
            ^{\frac{1}{t}}, \;\; i = 0, \cdots, l, \; t > 0
        \end{align}
        \noindent     in particular, if  $\rho$ is a faithful irrep of a group without a center, then  we have an improved "centralizer bound" ((cf. e.g. \cite{GLT1}, Introduction):
        \begin{align}
            \gamma_i \; \leq \;
            \left(\frac{   |G| - n^{2}  }  {
                |C_0| + \cdots + |C_i|  }  \right)
            ^{\frac{1}{2}}, \;\; i = 0, \cdots, l  \nonumber
        \end{align}
        \noindent     and
        \begin{align}
            \sqrt{ \frac{|G| - n^2  }{|G| -|G_0| -1} } \; \leq \;  \gamma \; \leq \; \sqrt{ \frac{|G| - n^2  }{|C_0|} }
        \end{align}
    \end{corollary}
    \noindent Our next corollary relates character values to Kronecker coefficients
    \begin{corollary}(cf. Lemma 1 (ii)),
        \begin{align}
            \gamma_i \; \leq \;
            \left(\frac{  \sum_{j=1}^k a_j^{2} - n^{4} |K|/|G| }  {
                ( |C_0| + \cdots + |C_i|  )/|G|}  \right)
            ^{\frac{1}{4}}, \;\; i = 1, \cdots, l    \nonumber
        \end{align}
        \noindent and in particular (cf. Remark 1.1(c))
        \begin{align}
            \sqrt[^4]{ \frac{|G| \cdot
                    \sum_{i=1}^k a_i^{2}
                    - |K|n^4  }{|G| - |G_0| - |K|}  }  \; \leq \; \gamma \;  \leq \; \sqrt[^4]{ \frac{|G| \cdot
                    \sum_{i=1}^k a_i^{2}
                    - |K|n^4  }{|C_0|}  }
        \end{align}

    \end{corollary}
    \begin{example}
        For a  $3$-dimensional irrep of the alternating group $A_5$  we have (cf. \cite{Serre2})
        $$
        \chi(1) =  n = 3, \; |G| = 60, \; |K| = 1, \; |C_0| = 12, \; |G_0| = 20 \nonumber, \;
        \gamma = (1 + \sqrt{5})/2, \; \norm{\chi\chi}^2 = 3 \nonumber
        $$
        Here the interval  (1.7) is
        $ \approx [1.14, \; 2.06 ] $
        and the interval (1.8) is $ \approx [1.26, \;  1.69]  $.

        \noindent The five-dimensional irrep of $A_5$ is an interesting corner case (ibid.) where
        \begin{align}
            \chi(1) =  n = 5, \; |G| = 60, \; |K| = 1, \; |C_0| = 35, \; , \; |G_0| = 24, \;
            \gamma = 1, \; \norm{\chi\chi}^2 = 11 \nonumber
        \end{align}
        \noindent and all bounds in (1.7) and (1.8) are equal to $1$ exactly
    \end{example}

    \section{Book-keeping, definitions and lemmas}
    Let  $  \mathfrak{c}_j, \; j=2, \cdots, k $ be all  conjugate classes of $ G $ in $ G \setminus K $. By (1.4)
    $$  a_{i,t}/n^{t} = \; < \chi_i, \; |\chi|^{t}/n^{t}   > $$
    and therefore
    \begin{align}
        a_{i,t}/n^{t} = (1/|G|)\left(n_i|K|  +  \frac{1}{n^{t}}\sum_{j=2}^k |\mathfrak{c}_j| \cdot \chi_i(\mathfrak{c}_j) \cdot |\chi(\mathfrak{c}_j)|^{t} \right), \; i = 1, \cdots, k
    \end{align}
    \noindent The sets  $ C_q , \; q = 0, \cdots,  l$  form a partition of $ G \setminus K $ ($C_l = G_0$) and in turn each of these sets is  a disjoint union of one or more conjugate classes $ \mathfrak{c}_j $. Take one of the irreducible characters $\chi_i $ from the list (1.3) (case $ \chi_i \equiv \chi $ is not excluded, cf. Remark 1.2) and define "incidence numbers" that depend only on a pair of irreducible characters $ \chi, \chi_i $
    \begin{align}
        \iota_{i, q} = \sum_{ \mathfrak{c} \subset C_q } |\mathfrak{c}| \cdot \chi_i(\mathfrak{c}) = \sum_{g \in C_q} \chi_i(g), \;\; i = 1, \cdots , k; \; q = 0, \cdots, l
    \end{align}
    where sum on the left is over all the conjugate classes $ \mathfrak{c} $ that are contained in $ C_q $.
    \noindent The following obvious lemma illustrates the definition (cf. (2.1-2))
    \begin{lemma}
        For $ i = 1, \cdots, k $ and any real $ t \geq 0 $
        \begin{align}
            a_{i,t}/n^{t} - (|K|/|G|)n_i  \; = \;
            \frac{1}{|G|} \sum_{q=0}^{l} \iota_{i,q} \left(\frac{\gamma_q}{n}\right) ^{t}
        \end{align}
    \end{lemma}

    \begin{remark}\
        \begin{enumerate}
            \item [(a)] Setting
            $ \iota_{i,K} = \sum_{g \in K }\chi_i(g) =  |K| \cdot n_i   $ for  $ i = 1, \cdots, k $ we note that the term   $ (1/|G|)\cdot|K|\cdot n_i$  conforms to the pattern of the right hind side of (2.3)  (cf. Remark 5.4)
            \item[(b)]
            The summation in (2.3) is actually up-to $ l -1 $ as $ \gamma_l = 0 $ (cf. (1.1))
        \end{enumerate}
    \end{remark}

    \begin{lemma}
        Numbers $  \iota_{i,q} $ and  $ a_{i,t} $  are real for all $  i = 1, \cdots , k; \; q = 0, \cdots, l-1 $ and $ t \geq 0$.   Numbers $ a_{i,t} $ are not necessarily non-negative (cf. Example 2.1). However, all numbers $ a_{i,t} $ are  positive  for sufficiently large $t$.
    \end{lemma}

    \noindent Proof. The numbers  $ a_{i,t} $  are integer when $t$ is an even integer. Hence, by (2.3)  $$ \sum_{q=0}^{l-1} \Im(\iota_{i,q}) \left(\frac{\gamma_q}{n}\right) ^{t} = 0 $$ for all even integers $t$, and therefore  $ \Im(\iota_{i,q}) = 0 $ for all
    $ q < l $ and all $i, \;i = 1, \cdots , k  $.  Further, it follows from (2.3) (cf. \cite{COUL}) that $\lim_{t \rightarrow \infty} a_{i,t}/n^t =  (|K|/|G|)n_i, \; i = 1, \cdots , k $ and therefore $ a_{i,t} > 0 $ for large $t$ (see an example below)
    \begin{example}
        Let  $\chi_3$ be  one of the  $3$-dimensional irreducible characters of $A_5$ (cf. Example 1.1).  By Lemma 2.1, the multiplicity $m(t)$ (see Figure 1) of the other $3$-dimensional irreducible character $ \chi'_3 $   of $A_5$ in virtual representation  $ (|\chi_3|/3)^t $ is
        \begin{align}
            m(t) =     3/60 - 3^{-t}\left[(15/60) - (12/60) (ab^t -ba^t)\right] \nonumber
        \end{align}
        where $ a = (\sqrt{5} + 1)/2, \; b = (\sqrt{5} - 1 )/2 $. It is easy to see that $m(t) < 0$   when $ 0 <  t < 2 $.

        The multiplicity of the $5$-dimensional irrep of $A_5$ in  $ (|\chi|/3)^t $ is
        $ 1/12 + (1/4)3^{-t} $ and
        we see  that the leading term $ (\gamma_0/n)^t $  is not necessarily present in  (2.3) (cf. Theorem 4.1).
    \end{example}
    \begin{figure}[!ht]
        \centering
        \includegraphics[ height=.5\linewidth]{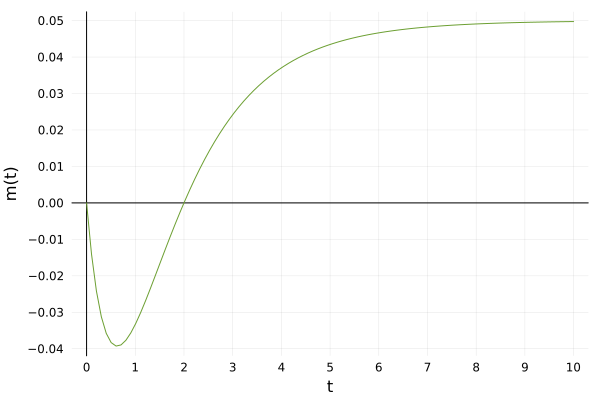}

        \caption{ Multiplicity $m(t)$ of $\chi'_3$ in $ (|\chi_3|/\chi_3(1))^t, \; t \geq 0 $ (cf. Example 2.1).  The graph is produced by  Julia Plots package (cf. \cite{JuliaPlots}). It seems that the most interesting part of the graph is located within the interval $0 \leq t \leq 2$. The effect  of $\chi'_3$ on  $ (|\chi_3|/\chi_3(1))^t $ is zero at both ends of the interval and is negative in its interior. The  effect of  $\chi'_3$  on the "tensor power"   $ (|\chi_3|/\chi_3(1))^t $ becomes  positive when $ t $ exceeds $2$}

    \end{figure}

    \begin{lemma}\
        \begin{enumerate}
            \item[(a)] If   the set $C_q$ is a conjugate class in $ G $ for some fixed $ q :  0 \leq q \leq l $ then \linebreak  $ \iota_{i, q} = |C_q| \chi_i(C_q) $ for all $ i = 1, \cdots, k $ (cf. Remark 1.1(b))
            \item[(b)] $  \iota_{1, q} = |C_q|, \; q = 0, \cdots, l $
            \item[(c)] For any fixed $i, \; 1 \leq i \leq k $ there
            is $ q, \; 0 \leq q \leq l$ such that $  \iota_{i, q} \neq 0 $
        \end{enumerate}
    \end{lemma}
    \noindent Statements (a) and (b) are obvious and the statement (c) is obvious  for $i = 1$.
    If $  i > 1 $ then
    $$ 0 = |G| \mathbb{E}_G(\chi_i) = |K|n_i + \sum_{q=1}^l   \iota_{i, q}  $$
\noindent by (2.2) and that proves (c)
    \begin{remark}
        The statement  (c) of Lemma 2.3  is not particularly useful.  It is unclear whether
        the numbers $  \iota_{i, q}, \; q = 0 \, \cdots, l-1 $ could all be  zero. Such cases, however, seem to be highly unlikely. The next lemma provides some obvious sufficient conditions
        for the right hand side of (2.3) to be nontrivial
    \end{remark}
    \begin{lemma}
        For a fixed $i, \; 1 \leq i \leq k $ there
        is $ q, \; 0 \leq q < l$ such that $  \iota_{i, q} \neq 0 $
        at least in any of the following cases
        \begin{enumerate}
            \item[(i)] $ n^2 n_i $ is not divisible by $ |G/K| $
            \item[(ii)] there is $ t \geq 0 $  such that $ a_{i,t}  \neq  (|K|/|G|)n^{t}n_i$, e.g. $ a_{i,t} \leq  0 $
            \item[(iii)] all sets $ C_q, \; q =0, \cdots, l-1$ are conjugate classes of $ G $ (in this case the irrep $\rho$ is necessarily real, cf. Remark 1.1 (b))
            \item[(iv)] $ \chi_i = \chi $ if $ \chi $ belongs to  the list (1.3) (cf, Remark 1.2)
        \end{enumerate}
    \end{lemma}
    \noindent To verify  the statement (i), note that  if the right hand side of (2.3) is zero then for the integer Kronecker coefficient  $a_i \equiv a_{i,2} $ we have $ a_i = (|K|/|G)|n^2 n_i $.
    Statements (ii)
    and (iii) are obvious.
    The  statement (iv) can be proved in the same way as Lemma 2.3 (c)
    \begin{conjecture}
        Under assumptions of Remark 1.1 (a) (i.e. $l>0$ in (1.1))  at least one of the incidence coefficients $  \iota_{i, q}, \; q = 0, \; \cdots, l-1 $ is non-zero  for any character $\chi_i$ on the list  (1.3).
    \end{conjecture}

    \noindent Formula (2.3) is much simpler for $ i =1 $ (cf.
    Lemma 2.3 (b))
    \begin{lemma} (Cf. \cite{He} and Lemma 2.3 (b)). For any real $ t \geq 0 $
        \begin{align}
            a_{1,t}/n^{t} - |K|/|G| \; = \;
            \frac{1}{|G|} \sum_{i=0}^{l} |C_i|\left(\frac{\gamma_i}{n}\right) ^{t}
        \end{align}
    \end{lemma}

    \section{Proof of Theorem 1.1 and Corollary 1.2}
    \noindent Fix $i, \; 0 \leq i \leq l - 1 $. By Lemma 2.5
    \begin{align}
        a_{1,t}/n^{t} - |K|/|G| \; \geq \;
        \frac{1}{|G|} \left(\frac{\gamma_i}{n}\right) ^{t}\sum_{j=0}^{i} |C_j|, \;\; i = 0, \cdots, l  \nonumber
    \end{align}
    \noindent and we get the upper bound in (1.5)  as well as  all upper bounds (1.6) of Corollary 1.2.
    Further, rewrite (2.4) as
    \begin{align}
        \frac {a_{1,t}/n^{t} - |K|/|G|} {|C_0|/|G|} \; = \;
        \left(\frac{\gamma}{n} \right) ^{t} \left(1 + \sum_{i=1}^{l} \frac{|C_i|}{|C_0|}\left(\frac{\gamma_i}{\gamma}\right) ^{t} \right)
    \end{align}
    \noindent and note that
    it follows from (3.1) that
    \begin{align}
        \frac {a_{1,t}/n^{t} - |K|/|G|} {|C_0|/|G|}
        \; \leq \;
        \left( \frac{\gamma}{n} \right) ^{t}
        \left(
        \frac{|C_0| + \cdots + |C_{l-1}|}{|C_0|}
        \right)             \nonumber
    \end{align}
    \noindent  The lower bound in statement (a) of Theorem 1.1 becomes now obvious if one recalls definition (1.1) and Remark 1.1 (c). Finally, (3.1) implies that
    \begin{align}
        \lim_{ t \rightarrow \infty } \left( a_{1,t}/n^{t} - |K|/|G| \right)^{1/t} \; = \;
        \frac{\gamma}{n}     \nonumber
    \end{align}
    which is equivalent to the statement (b) of Theorem 1.1. The statement (d) directly follows from  (3.1) as well

    \section{Multiplicities of irreps  in $(|\chi|/n)^t $ }
    We summarize some of the points of the preceding discussion by  self explanatory
    \begin{theorem} (cf. \cite{He}, \cite{COUL}).
        Fix $ i, \; 1 \leq i \leq k $ and  let   $ \iota_i = \iota_{i,q_i}$   be the first non-zero coefficient among $ \iota_{i,q}, \; q = 0, \cdots l - 1 $  (if it exists, cf. Example 2.1 and Remark 2.2). Let $ \hat{\gamma}_i = \gamma_{q_i} $ (cf. (1.1)) be the corresponding absolute value of the character $\chi$.
        \begin{enumerate}

            \item[(i)]  If all $ \iota_{i,q}, \; q = 0, \cdots l - 1 $ are zero (see Remark 2.2) then
            $  a_{i,t}/n^{t} = (|K|/|G|)n_i $  for any $ t \geq 0 $ and, in any case  (cf. Lemma 2.1 and \cite{COUL})
            $$
            \lim_{ t \rightarrow \infty } a_{i,t}/n^{t} = (|K|/|G|)n_i
            $$

            \item[(ii)] (cf. Lemma 2.1). $ a_{i,t}/n^{t}$ approaches $(|K|/|G|)n_i$ from below if $ \iota_i > 0 $ and from above if $ \iota_i < 0 $.
            In both cases $$
            \lim_{ t \rightarrow \infty } |\; a_{i,t}/n^{t} - (|K|/|G|)n_i\;|^{1/t} = \hat{\gamma}_i/n
            $$
            \item[(iii)] (cf. Lemma 2.1).
            $ ( \hat{\gamma}_i/n )^{t} $ and
            $  \iota_{i}^{-1} [a_{i,t}/n^{t} - (|K|/|G|)n_i] $
            are asymptotically equal as functions of $ t $
        \end{enumerate}
    \end{theorem}
    \begin{corollary}(cf. \cite{COUL}, \cite{He} and Lemma 2.2)
        $$ \lim_{t \rightarrow \infty}  (a_{i,t})^{1/t}  = n, \; i = 1, \cdots , k  $$
    \end{corollary}
    \noindent
    The statement (i) of Theorem 4.1, as well as Corollary 4.1 were obtained in a more general settings in \cite{COUL}.
    We reproduce here some  arguments from \cite{COUL}  in a slightly different context.
    Let   $ \delta_K $ denote the characteristic function  of the subset $ K $ of $G $

    \begin{proposition} (\cite{COUL})
        \begin{enumerate}
            \item[(i)]  $ \lim_{t\rightarrow \infty}  |\chi|^{t}/n^{t} =  \delta_K  $
            \item[(ii)]
            $ \delta_K = (|K|/|G|)\sum_{i=1}^{k} n_i \chi_i $
            \item[(iii)]
            $ \lim_{t \rightarrow \infty}  (a_{i,t}/n^{t})  = (|K|/|G|) n_i, \; i = 1, \cdots , k $
            \item[(iv)]
            $ \lim_{t \rightarrow \infty}  (a_{i,t})^{1/t}  = n, \; i = 1, \cdots , k  $
            \item[(v)] Let $ c_t(\rho) = \sum_{i=1}^k a_{i,t}/n^t $. Then the limit  $ c(\rho) = \lim_{t \rightarrow \infty}c_t(\rho) $ exists,
            \begin{align}
                c(\rho)= (|K|/|G|)\sum_{i=1}^k n_i
            \end{align}
            and
            \begin{align}
                c(\rho)  \; \leq \;  \sqrt{ \frac{k} {  |G/K|} }
            \end{align}

        \end{enumerate}

    \end{proposition}
    \noindent Indeed, the statement  (i) is obvious. To verify (ii), note that
  $
        \sum_{i=1}^{k} n_i \chi_i  = (|G|/|K|) \delta_K \nonumber
    $
    because $\sum_{i=1}^{k} n_i \chi_i$ is a lift (to $G$) of the regular character of $ G/K $.  Now
    (iii)  follows from the comparison of (ii) and (i) (cf. (1.4)) while (iv) and (v) (4.1) are obvious consequences of (iii).  Finally (v) (4.2) follows from (v) 4.1 by Cauchy-Schwarz inequality because $ \sum_{i=1}^k n_i^2 = |G/K| $  ( \cite{COUL})

    \begin{remark}
        For an exact irrep of a group  without a center  one can replace  $K$ with $\{1\} $ in all of the above
    \end{remark}

    \section{Character values and commuting probability}
    Let $\Gamma$ be  a finite group.
    By definition, the probability $c(g)$ of $ g \in \Gamma $ to be a commutator in $G$ (cf. e.g. \cite{Howe}, \cite{Serre2})  is

    \begin{align}
        c(g)  = \frac{| \{ [x,y] = g \; |  \; x,y \in \Gamma  \}| } { | \Gamma \times \Gamma |}, \;\; g \in \Gamma \nonumber
    \end{align}
    If  $ \{ \psi_1, \cdots, \psi_r \} $ is  the set of all pairwise inequivalent complex characters of $ \Gamma $ then by
    the well known formula of Frobenius (cf. \cite{Howe}, \cite{Serre2})

    \begin{align}
        c(  g )
        \; = \;  \frac{ \psi_1(g)/\psi(1) + \cdots +  \psi_r(g)/\psi_r(1)  }
        {|\Gamma|}, \;\; g \in \Gamma
    \end{align}

    \noindent For a central function $ \phi$ on $ \Gamma$ and an arbitrary $ g \in \Gamma $  define (cf. \cite{prob})
    a complex-valued random variable $ \xi_{\phi,g} $ on $ \Gamma \times \Gamma $ by the rule
    \begin{align}
        \xi_{\phi, g}(a,b) = \phi([a,b]^{-1}g)  \nonumber
    \end{align}

    \noindent The proof of the following lemma can be found in  \cite{prob}

    \begin{lemma} (cf. \cite{prob}, \cite{Serre2}).
        If $ \phi   $ is an irreducible character of $ \Gamma $ then
        \begin{align}
            \mathbb{E}_{\Gamma \times \Gamma} (  \xi_{\phi,g} ) \; = \; \frac{1}{(\dim \phi)^2}  \phi(g)
        \end{align}
    \end{lemma}

    \begin{lemma}  (cf. \cite{prob}).  If $ \phi   $ is an irreducible character of $ \Gamma $ then for any $ t > 0 $
        $$
        \mathbb{E}_{\Gamma \times \Gamma} ( | \xi_{\phi, g}|^t )/ n^t \; \geq \; c(g)
        $$
    \end{lemma}
    \noindent Proof. If $ g $ is a commutator $  | \xi_{\phi, g}|  = n  $  . Hence, by  Markov inequality (cf. \cite{B}) we get
    $$  c(g) = \textnormal{Pr}_{\Gamma \times \Gamma}( g \textnormal{ is a commutator} ) \leq \textnormal{Pr}_{\Gamma \times \Gamma}( |\xi_{\phi,g}| \geq n ) \; \leq \; \mathbb{E}_{\Gamma \times \Gamma} ( | \xi_{\phi, g}|^t )/ n^t $$
    \noindent Returning to our fixed irrep $\rho$ (1.0) we will write
    $ \xi_{\chi, g } \equiv \xi_{g}$
    assuming from now on that  $\rho$   is exact and  that the center of $ G $ is trivial (cf. Remark 1.2).
    \begin{remark}
        We thus assume that $ K = \{ 1 \} $ (cf. Remark 4.1). The assumption allows to simplify some formulas. Strictly speaking, this condition is not necessary since one can pass to the lifts of irreps of $ G /K$ as explained in Introduction
    \end{remark}
    \noindent The immediate consequence of these assumptions is

    \begin{lemma}
        $ c(g) =  \textnormal{Pr}_{G\times G}( |\xi_{g}| = n )$ for any $ g \in G $.   In other words,
        $ \xi_{g} (a,b) = n $  if and only if $ g = [a,b] $.

    \end{lemma}
    \noindent It turns out that a formula  similar to (5.1) still holds for  probabilities $ \textnormal{Pr}_{G\times G}( |\xi_{g}| = \gamma_q)$ with any $ q = 0, \cdots, l  $.  We need, however, to outline some preliminary details   before   making a precise statement.  In particular, we would like to estimate the expectation $\mathbb{E}_{G \times G} ( | \xi_g|^t )$
    (cf. \cite{prob}).

\noindent     First of all, from Lemma 5.1 and definition (1.4) we get
    \begin{lemma} (cf. \cite{prob}). For any $ t \geq 0 $
        \begin{align}
            \mathbb{E}_{G \times G} ( | \xi_g|^t ) \; = \;
            \sum_{i=1}^k (a_{i,t}/n_i^2) \chi_i(g)
        \end{align}
    \end{lemma}

    \noindent On the other hand, we have

    \begin{lemma}\
        \begin{enumerate}
            \item[(i)]
            \begin{align}
                \mathbb{E}_{G \times G} ( | \xi_g|^t )/ n^t \; - \;  c(g) \; = \;
                \frac{1}{|G|^2} \sum_{a,b \in G; \;[a ,b] \neq g } \frac{ |\chi([a,b]^{-1}g)|^t}{n^t}
            \end{align}
            \item[(ii)] There are non-negative integers $ b_q, \; q = 0, \cdots , l -1 $ that do not depend on $ t $ such that
            \begin{align}
                \mathbb{E}_{G \times G} ( | \xi_g|^t )/ n^t - c(g) =
                \frac{1}{|G|^2} \sum_{q=0}^{l-1} b_q \left(\frac{\gamma_q}{n} \right)^t
            \end{align}
            \item[(iii)] $ \textnormal{Pr}_{G\times G}( |\xi_{g}| \; = \; \gamma_q)  \; = \; b_q / |G|^2 , \; q = 0, \cdots, l-1  $
        \end{enumerate}
    \end{lemma}
    \noindent To prove the statement (i) write
    \begin{align}
        \mathbb{E}_{G \times G} ( | \xi_g|^t ) \; = \;
        \frac{1}{|G|^2} \left( \; | \{ (a,b) \in G \times G , \; [a,b] = g \} | \cdot n^t  \;\; +   \sum_{a,b \in G; \;[a ,b] \neq g }  |\chi([a,b]^{-1}g)|^t \right) \nonumber
    \end{align}
    and take Lemma 5.3 into account. To verify the statement (ii) and (iii) just collect equal terms on the right hand side of (5.4)
    \begin{remark}
        The right hand side of (5.4) is non-zero if every element of $ G $ s a commutator (e.g. if $ G $ is a simple group - cf. \cite{Ore}). Indeed, let
        $$ H =
        g \cdot \{ [a,b]^{-1}, \; a,b \in G \; | \;  [a,b] \neq g   \}$$
        If every element of $ G $ is a commutator then $  H = G \setminus \{1\} $ and there is $ h \in H $ such that $ \chi(h) \neq 0 $. One could ask if the reverse  of this statement is also true.

        \begin{conjecture}
            Statements
            "$ \mathbb{E}_{G \times G} ( | \xi_g| )/ n \; \neq \; c(g) $ for every $ g \in G $" and
            "every element of $ G $ is a commutator" are equivalent
        \end{conjecture}

    \end{remark}
    \begin{remark}
        It is easy to see that the  Conjecture 5.1 is not true for the excluded case of (monomial) irreps of nilpotent groups (cf. Remark 1.1(a))
    \end{remark}
    \noindent
    Plugging (2.3) into (5.3) and rearranging terms (recall that here $|K| = 1$) we get
    \begin{align}
        \mathbb{E}_{G \times G} ( | \xi_g|^t )/ n^t \; = \;
        \frac{1}{|G|} \sum_{i=1}^k \frac{\chi_i(g)}{n_i}
        \;  + \;   \frac{1}{|G|} \sum_{q=0}^{l-1}  \left( \sum_{i=1}^k \frac{\chi_i(g) }{n_i^2} \iota_{i,q} \right) \left(\frac{\gamma_q}{n}\right)^t
    \end{align}
    \noindent By Frobenius formula (5.1) the first summand in (5.6) is precisely $ c(g)$ while coefficients in front of $( \gamma_q/n)^t $ do not depend on $ t $. Hence, by comparing (5.6) with
    Lemma 5.5 (ii)-(iii) and plugging in the values of $ \iota_{i,q} $ from (2.2) we obtain the formula
    \begin{align}
        \textnormal{Pr}_{G\times G}( |\xi_{g}| \; = \; \gamma_q) \; = \;
        \frac{1}{|G|}   \sum_{i=1}^k \frac{\chi_i(g) }{n_i^2} \iota_{i,q} \; = \;  \frac{1}{|G|}  \sum_{h \in C_q }  \sum_{i=1}^k    \frac{\chi_i(h) \chi_i(g) }{n_i^2}, \; g \in G
    \end{align}
    \noindent for all $  q = 0, \cdots, l-1 $.
    As suggested by the right hand side of (5.7) set
    \begin{align}
      c_q(g) = \frac{1}{|G|}  \sum_{h \in C_q }  \sum_{i=1}^k    \frac{\chi_i(h) \chi_i(g) }{n_i^2}, \; g \in G , \; q = 0, \cdots, l      \tag{5.7'}
    \end{align}
    \begin{remark}
        Applying (5.7) to the level set $  K  = \{1\} $ (cf. Remark 2.1 and Lemma 5.3) we recover Frobenius formula (5.1)
        $$ \!\!\!\!   \textnormal{Pr}_{G\times G}( |\xi_{g}| \; = \; n )  \; = \;
        \frac{1}{|G|}    \sum_{i=1}^k \frac{\chi_i(g) }{n_i^2} \iota_{i,\{1\}} \; = \;  \frac{1}{|G|}  \sum_{h \in \{1\} }  \sum_{i=1}^k    \frac{\chi_i(h) \chi_i(g)   }{n_i^2}   \; = \;         \frac{1}{|G|} \sum_{i=1}^k    \frac{ \chi_i(g)   }{n_i} \; = \;  c(g)                               $$
        \noindent
    \end{remark}

    \begin{theorem} For an exact irrep $ \rho $ (1.0) of a group $ G $ without a center and for any $ t \geq 0, \;  g \in G $
        \begin{enumerate}

            \item[(i)]
            $    \mathbb{E}_{G \times G} ( | \xi_g|^t )/ n^t \; \geq \; c(g) $

            \item[(ii)]
            $   \lim_{t \rightarrow \infty }     \mathbb{E}_{G \times G} ( | \xi_g|^t )/ n^t \; = \; c(g)
            $
            \item[(iii)]
            There is $ q, \; 0 \leq q \leq l $ (that depends on $ g $) such that
            $$
            \lim_{t \rightarrow \infty }  \left(   \mathbb{E}_{G \times G} ( | \xi_g|^t )/ n^t - c(g) \right)^{1/t} = \gamma_q/n
            $$
            \noindent and   functions $ \mathbb{E}_{G \times G} ( | \xi_g|^t )/ n^t - c(g)$ and $ c_q(g)  (\gamma_q/n)^t $ are asymptotically equal
            \item[(iv)] for all $ q = 0, \cdots, l $ (recall that $ \gamma_l = 0, \;  C_l = G_0 = |\chi|^{-1}(0)$) and $ g \in G $
        \begin{align}
\!\!\!\! \!\!\!\!    \!\!\!\!   \!\!\!\!    \!\!\!\!   \!\!\!\!        c_q(g) \; = \; \textnormal{Pr}_{G\times G}( |\xi_{g}| \; = \; \gamma_q)  \; = \;
           \frac{1}{|G|}   \sum_{i=1}^k \frac{\chi_i(g) }{n_i^2} \iota_{i,q} \; = \;  \frac{1}{|G|}  \sum_{h \in C_q }  \sum_{i=1}^k    \frac{\chi_i(h) \chi_i(g) }{n_i^2}   \nonumber
        \end{align}

        \end{enumerate}
    \end{theorem}
    \noindent Proof. The statement (i) follows from (5.5) (or in more general setting from Lemma 5.2). If there is $ b_q \neq 0, \;  0 \leq q < l $ in (5.5) then  the
    term $ (b_q/|G|^2)(\gamma_q/n)^t $ dominates the sum $ (1/|G|^2) \sum_{j=q}^{l} b_j (\gamma_j/n)^t $ when $ t $ is large enough. In other words,  it follows from (5.5) that unless  $ b_j =0  $  for all $ j = 0, \cdots, l $ there is an index $r, \; 0 \leq r < l $   such that $ b_{r} > 0  $ and $$   \mathbb{E}_{G \times G} ( | \xi_g|^t )/ n^t - c(g) = (b_{r}/|G|^2)  (\gamma_r/n)^t + c (\gamma_{r+1}/n)^t   $$ where  $ c $ is bounden when $t \rightarrow \infty$. That proves (ii) and (iii). For $ q <l $ the statement (iv) was already established (cf. (5.7)). In the remaining case of  $ q = l$ we argue as follows. By (5.7'), the sum of all numbers $ c(g) $ and  $ c_q(g), \; q = 0, \cdots, l, \; g \in G $ is equal to
    \begin{align}
        \sum_{i=1}^k \frac{n_i\chi_i(g)}{n_i^2} \;  + \; \frac{1}{|G|}  \sum_{h \in G, h \neq 1 }  \sum_{i=1}^k    \frac{\chi_i(h) \chi_i(g) }{n_i^2}  \; =  \;
        \frac{1}{|G|}   \sum_{i=1}^k \left( \sum_{h \in G } \chi_i(h) \right) \frac{ \chi_i(g) }{n_i^2}   \nonumber
    \end{align}
    \noindent The interior sum on the right hand side of this equation is equal to zero for all $i$ except for $ i = 1 $.  When $i=1$  the interior sum equals  $ |G| $ while
    $ \chi_i(g)/n_i^{2} = 1 $ for any $ g \in G $. The net result is that $     c(g) + \sum_{q=0}^l c_q(g)  =      1 $. On the other hand,  $ c(g) =  \textnormal{Pr}_{G\times G}( |\xi_{g}| = n )$ and $ c_q(g) =    \textnormal{Pr}_{G\times G}( |\xi_{g}| \; = \; \gamma_q) $ for all  $ g \in G , \; q = 0, \cdots, l-1 $  (Lemma 5.3 and (5.7)). Hence we have to conclude that (iv) is valid for $ q = l $ as well.
\newline\newline\noindent
By setting $ g =  1 $ in Theorem 5.1 (iv) we get
\begin{corollary} For all $ q = 0, \cdots,  l $
    \begin{align}
    c_q(1) = \textnormal{Pr}_{(a,b) \in G\times G}( |\chi([a,b])| = \gamma_q ) = \frac{1}{|G|}  \sum_{h \in C_q }  \sum_{i=1}^k    \frac{\chi_i(h) }{n_i}  =  \sum_{\mathfrak{c} \subset C_q } \frac{\mathfrak{|c|}}{|G|} \sum_{i=1}^k    \frac{\chi_i(\mathfrak{c}) }{n_i}
    \end{align}

\noindent where exterior sum is over all conjugate classes contained in $C_q$.
In particular,
 replacing here the level set $  C_q $ with the trivial level set $ \{1\} $ we get the "Frobenius formula" $k/|G|$ for the probability of a commuting pair of group elements (cf. Remark 5.4)
 \end{corollary}
 \begin{remark}
 Corollary 5.1 has an obvious   "physical meaning". The left hand side
   of (5.8) is the probability of a commutator to belong to a level set of the function $ |\chi| $. On the other hand, by  Frobenius formula (5.1), the sum $ \sum_{i=1}^k    \chi_i(h)/n_i $ is the probability of an element of the same  level set to be a commutator
 \end{remark}
    \begin{example}
        Level sets of the five-dimensional irreducible character $\chi_5 $  of $ A_5 $ ((cf. \cite{Serre2} and Example 2.1) are:
        $$     K = |\chi_5|^{-1}(n) = \{1\}, \;  |K| = 1;  \;
        C_0 = |\chi_5|^{-1}(1), \;  |C_0| = 35;  \;
        C_1 = |\chi_5|^{-1}(0) = G_0, \;  |C_1| = 24
        $$
        \noindent Using the recipe provided by Corollary 5.1 we get
        $$
        \textnormal{Pr}( |\chi([a,b])| = 5 ) = 5/60; \;
        \textnormal{Pr}( |\chi([a,b])| = 1 ) = 29/60; \;
        \textnormal{Pr}( |\chi([a,b])| = 0 ) = 26/60
        $$
    \end{example}
    \begin{cremark}
        Alongside \textbf{the partition} of a finite group  by conjugate classes there are also partitions
        by level sets of various central functions. Results presented above seem to indicate that these coarser partitions could be useful. Another  application of  partitions of this kind can be found in \cite{Id}
    \end{cremark}
    
\end{document}